\documentclass[letterpaper]{article}

\newcommand{\myauthor}{Benjamin Antieau}
\newcommand{\mytitle}{Twisted derived equivalences for affine schemes}
\newcommand{\pdftitle}{\mytitle}
\author{\myauthor}

\title{\mytitle}

\usepackage[pdfstartview=FitH,
            pdfauthor={\myauthor},
            pdftitle={\pdftitle},
            colorlinks,
            linkcolor=reference,
            citecolor=citation,
            urlcolor=e-mail,
            ]{hyperref}

\usepackage{mathrsfs}
\usepackage[mathscr]{euscript}
\usepackage{amsmath}
\usepackage{amscd}
\usepackage{amsbsy}
\usepackage{amssymb}
\usepackage{microtype}
\usepackage{enumerate}
\usepackage[all,cmtip]{xy}
\usepackage{tikz}
\usetikzlibrary{matrix,arrows}
\usepackage{dsfont}
\usepackage[abbrev,lite]{amsrefs}
\usepackage{amsthm}

\usepackage{color}
\definecolor{todo}{rgb}{1,0,0}
\definecolor{conditional}{rgb}{0,1,0}
\definecolor{e-mail}{rgb}{0,.40,.80}
\definecolor{reference}{rgb}{.20,.60,.22}
\definecolor{mrnumber}{rgb}{.80,.40,0}
\definecolor{citation}{rgb}{0,.40,.80}

\newcommand{\perf}{\mathrm{perf}}

\newcommand{\we}{\simeq}
\newcommand{\iso}{\cong}

\newcommand{\Gm}{\mathds{G}_{m}}

\newcommand{\caldararu}{C\u{a}ld\u{a}raru}

\newcommand{\tors}{\mathrm{tors}}

\newcommand{\qc}{\mathrm{qc}}

\newcommand{\Mod}{\mathrm{Mod}}
\newcommand{\Perf}{\mathrm{Perf}}

\newcommand{\QC}{\mathrm{QC}}

\DeclareMathOperator{\Spec}{Spec}

\DeclareMathOperator{\Hoh}{H}

\DeclareMathOperator{\End}{End}

\newcommand{\StQC}{\mathscr{QC}}

\newcommand{\et}{\mathrm{\acute{e}t}}

\DeclareMathOperator{\Br}{Br}
\DeclareMathOperator{\dBr}{dBr}

\newcommand{\Drm}{\mathrm{D}}

\newcommand{\Fscr}{\mathscr{F}}
\newcommand{\Oscr}{\mathscr{O}}

\newcommand{\Ascr}{\mathscr{A}}

\newcommand{\Xscr}{\mathscr{X}}

\newcommand{\CC}{\mathds{C}}

\newcommand{\ZZ}{\mathds{Z}}

\newcommand{\PP}{\mathds{P}}

\theoremstyle{plain}
\newtheorem{theorem}{Theorem}[section]

\newtheorem{corollary}[theorem]{Corollary}

\newtheoremstyle{named}{}{}{\itshape}{}{\bfseries}{.}{.5em}{#1 \thmnote{#3}}
\theoremstyle{named}

\theoremstyle{definition}

\newtheorem{question}[theorem]{Question}

\newtheorem{motto}[theorem]{Motto}
\newtheorem{remark}[theorem]{Remark}
\newtheorem{problem}[theorem]{Problem}

\let\oldmarginpar\marginpar
\renewcommand\marginpar[1]{\-\oldmarginpar[\raggedleft\footnotesize #1]%
{\raggedright\footnotesize #1}}

\usepackage{txfonts}

\begin{document}
\maketitle

\begin{abstract}
    \noindent
    We show how work of Rickard and To\"en completely resolves the question of when two
    twisted affine schemes are derived equivalent.

%
%

\end{abstract}

\section{Introduction}

The question of when $\Drm^b(X)$ is equivalent as a $k$-linear triangulated category to
$\Drm^b(Y)$ for two $k$-varieties $X$ and $Y$ has been extensively studied since Mukai
proved that $\Drm^b(\hat{A})\we\Drm^b(A)$ for an abelian variety $A$ and its dual
$\hat{A}$~\cite{mukai}.
Since in general $A$ and $\hat{A}$ are not isomorphic, derived equivalence is a weaker
condition than isomorphism. However, derived equivalence nevertheless does preserve a great deal
of information: derived equivalent varieties have the same dimension, the same algebraic
$K$-theory, and the same Hochschild homology.

The cohomological Brauer group of a scheme $X$ is
$\Br'(X)=\Hoh^2_{\et}(X,\Gm)_{\tors}$. When $X$ is quasi-compact, there is an inclusion
$\Br(X)\subseteq\Br'(X)$, where $\Br(X)$ denotes the Brauer group of $X$,
which classifies Brauer equivalence classes of Azumaya algebras on $X$. In many cases of
interest, $\Br(X)=\Br'(X)$.
Examples include all quasi-projective schemes over affine schemes~\cite{dejong}. The
Brauer group comes into play because in many problems on moduli of vector bundles, there is
an obstruction, living in the Brauer group of the coarse moduli space, to the existence of a
universal vector bundle. Another way to say this is that this class in the Brauer group is
the obstruction to the coarse moduli space being fine. At times one then obtains an equivalence
$\Drm^b(X)\we\Drm^b(Y,\beta)$, where $\Drm^b(Y,\beta)$ is the derived category of
$\beta$-twisted coherent sheaves. Particular cases of this arise in the study of K3 surfaces
for example.

The systematic study of when $\Drm^b(X,\alpha)\we\Drm^b(Y,\beta)$ began with
\caldararu's thesis~\cite{caldararu}. In this short note, we are interested in the following
two problems.

\begin{problem}\label{problem:A}
    Let $R$ be a commutative ring, 
    let $X$ and $Y$ be two locally noetherian $R$-schemes, and fix $\alpha\in\Br'(X)$ and $\beta\in\Br'(Y)$.
    Determine when there exists an $R$-linear equivalence of triangulated categories
    $\Drm^b(X,\alpha)\we\Drm^b(Y,\beta)$.
\end{problem}

\begin{problem}\label{problem:A'}
    Let $R$ be a commutative ring, 
    let $X$ and $Y$ be two quasi-compact and quasi-separated $R$-schemes, and fix $\alpha\in\Br'(X)$ and $\beta\in\Br'(Y)$.
    Determine when there exists an $R$-linear equivalence of triangulated categories
    $\Drm_{\perf}(X,\alpha)\we\Drm_{\perf}(Y,\beta)$.
\end{problem}

Here, $\Drm^b(X,\alpha)$ denotes the bounded derived category of
$\alpha$-twisted coherent sheaves, while $\Drm_{\perf}(X,\alpha)$ is the triangulated
category of perfect complexes of $\alpha$-twisted $\Oscr_X$-modules. When $X$ is regular and
noetherian, the existence locally of finite-length finitely generated locally free
resolutions implies that the natural map
$\Drm_{\perf}(X,\alpha)\rightarrow\Drm^b(X,\alpha)$ is an equivalence of $R$-linear triangulated
categories.

\begin{question}\label{question:eq}
    Are Problems~\ref{problem:A} and~\ref{problem:A'} equivalent for $X$ and $Y$ noetherian
    and quasi-separated?
\end{question}

The contents of our paper are as follows. In Section~\ref{sec:tdc}, we give some background
on twisted derived categories and equivalences between them.
Then, in Section~\ref{sec:affine}, the affine case of Problems~\ref{problem:A}
and~\ref{problem:A'} is solved completely, and we explain how work of Rickard shows that
these two problems \emph{are} equivalent for affine schemes. We do not claim that this
result is new, but rather that it is not as well-known as it should be.

\paragraph{Acknowledgments.}
We would like to thank the organizers of the AIM workshop ``Brauer groups and
obstruction problems'' for facilitating a stimulating week in February 2013.

\section{Twisted derived categories}\label{sec:tdc}

Let $X$ be a scheme, and take $\alpha\in\Hoh^2_{\et}(X,\Gm)$. Then, $\alpha$ is represented
by a $\Gm$-gerbe $\Xscr\rightarrow X$. There is a good notion of quasi-coherent sheaf on
$\Xscr$, or of coherent sheaf when $X$ is locally noetherian. An $\Oscr_{\Xscr}$-module
$\Fscr$ comes naturally with a left action of the sheaf $\mathds{G}_{m,X}$. But, there is a
second, inertial action, which can be described as saying that a section
$u\in\mathds{G}_{m,X}(U)$ over $U\rightarrow\Xscr$ acts on $\Fscr(U)$ via the
isomorphism $u^*\Fscr_U\rightarrow\Fscr_U$, which induces an isomorphism
$u^*:\Fscr(U)\rightarrow\Fscr(U)$. There is an associated left action of the inertial
action. An $\alpha$-twisted $\Oscr_X$-module is by definition an
$\Oscr_{\Xscr}$-module $\Fscr$ for which these two left actions agree. It is shown
in~\cite{lieblich-moduli}*{Proposition 2.1.3.3} that this agrees with the definition of
$\alpha$-twisted sheaf given by \caldararu.

If $X$ is a scheme and $\alpha\in\Hoh^2_{\et}(X,\Gm)$,
write $\Drm_{\perf}(X,\alpha)$ for the derived category of complexes of $\alpha$-twisted
sheaves that are \'etale locally quasi-isomorphic to finite-length complexes of vector
bundles. This is naturally a full subcategory of $\Drm_{\qc}(X,\alpha)$ of complexes of
$\alpha$-twisted sheaves with ($\alpha$-twisted) quasi-coherent cohomology sheaves. If $X$
is regular and noetherian, then $\Drm_{\perf}(X,\alpha)\we\Drm^b(X,\alpha)$, the bounded
derived category of $\alpha$-twisted coherent $\Oscr_X$-modules.

Let $A$ be an Azumaya algebra on $X$ with class $\alpha$. A complex of right $A$-modules $P$ (in the
abelian category $\Mod_{\Oscr_X}$) is perfect if there is an open affine cover $\{U_i\}_{i\in
I}$ of $X$ such that $P_{U_i}$ is quasi-isomorphic to a bounded complex of finitely
generated projective right $\Gamma(U_i,A)$-modules. The derived category of perfect complexes of
right $A$-modules will be denoted $\Drm_{\perf}(X,A)$. Then, as explained in~\cite{caldararu},
$\Drm_{\perf}(X,\alpha)\we\Drm_{\perf}(X,A)$. In the same way, there is a big derived
category of all all complexes of right $A$-modules with quasi-coherent cohomology sheaves $\Drm_{\qc}(X,A)$
and an equivalence $\Drm_{\qc}(X,\alpha)\we\Drm_{\qc}(X,A)$.

In the next section, we will need dg enhancements of these categories.
Write $\Perf(X,\alpha)$ and $\QC(X,\alpha)$ for dg enhancements of $\Drm_{\perf}(X,\alpha)$
and $\QC(X,\alpha)$, respectively. These are pretriangulated dg categories. The big dg
category $\QC(X,\alpha)$ is constructed for example in To\"en~\cite{toen-derived}. The small
dg category $\Perf(X,\alpha)$ can then be taken to be the dg category of compact objects in
$\QC(X,\alpha)$.

%
%
%

\section{Twisted derived equivalences over affine schemes}\label{sec:affine}

Many of us first learned of twisted derived categories from \caldararu's
thesis~\cite{caldararu} and the paper~\cite{caldararu-elliptic}. In that paper, \caldararu~cites a private communication from
Yekuteli giving the following theorem~\cite{caldararu}*{Theorem 6.2}.

\begin{theorem}
    Suppose that $R$ is a commutative local ring and that $A$ and $B$ are Azumaya
    $R$-algebras with classes $\alpha,\beta\in\Br(R)$. Then, the following are equivalent:
    \begin{enumerate}
        \item   $\alpha=\beta$ in $\Br(R)$;
        \item   $A$ and $B$ are derived Morita equivalent over $R$---that is, there is an
            $R$-linear equivalence of triangulated categories $\Drm^b(A)\we\Drm^b(B)$.
    \end{enumerate}
\end{theorem}

It is the purpose of our paper to advertise the fact that the condition that $R$ be local
is unnecessary.

\begin{theorem}\label{thm:thm}
    Suppose that $R$ is a commutative ring and that $\alpha,\beta\in\Br(R)$. Then, the
    following are equivalent:
    \begin{enumerate}
        \item   $\alpha=\beta$ in $\Br(R)$;
        \item   there is an $R$-linear equivalence of triangulated categories $\Drm_{\perf}(R,\alpha)\we\Drm_{\perf}(R,\beta)$.
    \end{enumerate}
    Moreover, if $R$ is noetherian, these are equivalent to:
    \begin{enumerate}
        \item[3.]   there is an $R$-linear equivalence of triangulated categories
            $\Drm^b(R,\alpha)\we\Drm^b(R,\beta)$.
    \end{enumerate}
    \begin{proof}
        Since $\alpha,\beta\in\Br(R)$, we can assume that $\alpha$ is represented by an
        Azumaya $R$-algebra $A$, and that $\beta$ is represented by an Azumaya $R$-algebra
        $B$. In this case,
        $\Drm_{\perf}(X,\alpha)\we\Drm_{\perf}(X,A)\we\Drm^b(\mathrm{proj}_A)$, where
        $\Drm^b(\mathrm{proj}_A)$ is the bounded derived category of finitely generated
        projective $A$-modules. The second equivalence follows because on an affine scheme,
        every perfect complex is quasi-isomorphic to a bounded complex of finitely generated
        projectives (see~\cite{thomason-trobaugh}*{Proposition~2.3.1(d)}), and this generalizes in a straightforward way to Azumaya algebras on an
        affine scheme.

        When $\alpha=\beta$, the Azumaya algebras $A$ and $B$ are Brauer equivalent. This
        means that there exist finitely generated projective $R$-modules $M$ and $N$ and an
        $R$-algebra isomorphism
        \begin{equation*}
            A\otimes_R\End_R(M)\iso B\otimes_R\End_R(N).
        \end{equation*}
        It follows from classical Morita theory that there is an equivalence
        $\Mod_A\we\Mod_B$ of abelian categories of right $A$ and $B$-modules. From this it
        follows immediately that $\Drm^b(\mathrm{proj}_A)\we\Drm^b(\mathrm{proj}_B)$. This
        proves that (1) implies (2).

        So, suppose that
        $\Drm_{\perf}(R,\alpha)\we\Drm_{\perf}(R,\beta)$, or in other words that
        $\Drm^b(\mathrm{proj}_A)\we\Drm^b(\mathrm{proj}_B)$. Rickard's
        theorem~\cite{rickard-morita}*{Theorem 6.4} as refined in~\cite{rickard-derived}
        implies that there is a tilting complex inducing an
        $R$-linear equivalence $\Drm^b(\mathrm{proj}_A)\we\Drm^b(\mathrm{proj}_B)$. (Rickard's
        theorem does not imply that this is the equivalence we began with, but it \emph{is}
        still $R$-linear.) The existence of
        the tilting complex implies that there is an equivalence of $R$-linear dg category enhancements
        $\Perf(R,\alpha)\we\Perf(R,\beta)$, which is then a derived Morita equivalence. That
        is, there is an equivalence of the ``big'' $R$-linear dg categories
        $\QC(R,\alpha)\we\QC(R,\beta)$. These are locally presentable dg categories with
        descent in the language of~\cite{toen-derived}. Now, the derived Brauer group of
        $R$, denoted $\dBr(R)$, classifies locally presentable dg categories with descent
        over $R$ that are \'etale locally equivalent to $\QC(R)$. Since $\Spec R$ is affine,
        the $R$-linear equivalence $\QC(R,\alpha)\we\QC(R,\beta)$ means that $\alpha$ and $\beta$
        define the same element of $\dBr(R)$ (see~\cite{toen-derived}*{Section 3}). But, $\dBr(R)\iso\Hoh^1_{\et}(\Spec
        R,\ZZ)\times\Hoh^2_{\et}(\Spec R,\Gm)$ by~\cite{toen-derived}*{Theorem 1.1}. Since $\Br(R)\subseteq\dBr(R)$, it follows
        that $\alpha=\beta$, and so (2) implies (1).

        Finally, the fact that (2) and (3) are equivalent follows
        from~\cite{rickard-morita}*{Propositions 8.1, 8.2}. This completes the proof.
    \end{proof}
\end{theorem}

\begin{remark}
    By~\cite{dejong}, $\Br(R)=\Br'(R)=\Hoh^2_{\et}(\Spec R,\Gm)_{\tors}$.
\end{remark}

We expand briefly on the philosophy of the proof.
Write $\StQC$ for the \'etale stack of locally presentable dg categories with dg category of
sections over $f:Y\rightarrow X$ the locally presentable dg category $\QC(Y)$, which is a dg
category enhancement of the derived category $\Drm_{\qc}(Y)$ of complexes of
$\Oscr_Y$-modules with quasi-coherent cohomology sheaves. The derived Brauer group $\dBr(X)$ of a scheme
classifies stacks of locally presentable dg categories that are \'etale locally equivalent
to $\StQC$ up to equivalence \emph{of stacks}.

\begin{motto}
    The Brauer group classifies Azumaya algebras $\Ascr$ up to derived Morita equivalence of stacks
    of dg categories of complexes of $\Ascr$-modules.
\end{motto}

For $\alpha\in\dBr(X)$, write $\StQC(\alpha)$ for the associated stack. For instance, if
$\alpha$ is the Brauer class of an Azumaya algebra $\Ascr$ over $X$
then the dg category of sections over $f:Y\rightarrow X$ of $\StQC(\alpha)$ is
$\QC(Y,f^*\Ascr)$, which is a dg category enhancement of
$\Drm_{\qc}(Y,\Ascr)\we\Drm_{\qc}(Y,\alpha)$. The key point in the proof of the theorem was
that over an affine scheme $\Spec R$, giving an equivalence of stacks
$\StQC(\alpha)\we\StQC(\beta)$ is equivalent to giving an $R$-linear equivalence of the global sections $\QC(R,\alpha)\we\QC(R,\beta)$.

On non-affine schemes, giving an equivalence of global sections is, not surprisingly, insufficient.
The following example is due to \caldararu~\cite{caldararu}*{Example 1.3.16}. Let $X$ be a smooth projective K3 surface over
the complex numbers given as a
double cover of $\PP^2$ branched along a smooth sextic curve. The involution $\phi$ of $X$ given by interchanging
the sheets of the cover has the property that $\phi^*\alpha=-\alpha$ for $\alpha\in\Br(X)$.
Clearly $\phi$ induces an equivalence
$\Drm^b(X,\alpha)\we\Drm^b(X,-\alpha)$. But, since $\Br(X)$ contains non-zero $p$-torsion
for every prime $p$, there is a class $\alpha\in\Br(X)$ such that $\alpha\neq-\alpha$. Thus,
the theorem fails in the non-affine case. The problem is that the equivalence does not
respect restriction to open subsets of $X$.

We now prove the conjecture suggested by \caldararu~after~\cite{caldararu-elliptic}*{Theorem
6.2}.

\begin{corollary}
    Suppose that $R$ and $S$ are commutative rings and that there is an equivalence of
    triangulated categories
    $\Drm_{\perf}(R,\alpha)\we\Drm_{\perf}(S,\beta)$ for $\alpha\in\Br(R)$ and
    $\beta\in\Br(S)$. Then, there exists a ring isomorphism $\phi:R\rightarrow S$ such that
    $\phi^*(\alpha)=\beta$ in $\Br(S)$.
    \begin{proof}
        Let $A$ be an Azumaya algebra with class $\alpha$ over $R$, and let $B$ an be Azumaya
        over $S$ with class $\beta$. Then, our hypotheses say that
        $\Drm^b(\mathrm{proj}_A)\we\Drm^b(\mathrm{proj}_B)$.
        By Rickard~\cite{rickard-morita}*{Proposition 9.2}, the centers of $A$ and $B$ are
        isomorphic.
        Thus, there is an isomorphism $\phi:R\rightarrow S$, and there are equivalences
        $\Drm_{\perf}(S,\beta)\we\Drm_{\perf}(R,\alpha)\we\Drm_{\perf}(S,\phi^*(\alpha))$.
        The composition induces a ring automorphism $\sigma:S\rightarrow S$. So, by composing on
        the $\phi$ on the right with $\sigma^{-1}$,
        we can assume that
        $\Drm_{\perf}(S,\beta)\we\Drm_{\perf}(S,\phi^*(\alpha))$ \emph{is} $S$-linear.
        The corollary follows now from the theorem.
    \end{proof}
\end{corollary}

We end by observing that the condition of $R$-linearity in Theorem~\ref{thm:thm} is necessary.

\begin{remark}\label{rem:caution}
    Consider the field $k=\CC(w,x,y,z)$ and the quaternion division algebras
    $(w,x)$ and $(y,z)$ over $k$. Then, these algebras are evidently derived Morita
    equivalent over $\CC$ (they are even isomorphic over $\CC$). However, $[(w,x)]\neq [(y,z)]$ in $\Br(k)$.
\end{remark}


\begin{bibdiv}
\begin{biblist}

%

    



%




















 





 
 


 



\bib{caldararu}{thesis}{
    author = {C\u{a}ld\u{a}raru, Andrei},
    title = {Derived categories of twisted sheaves on Calabi-Yau manifolds},
    note = {Ph.D. thesis, Cornell University (2000)},
    eprint = {http://www.math.wisc.edu/\textasciitilde andreic/},
}

\bib{caldararu-elliptic}{article}{
    author={C\u{a}ld\u{a}raru, Andrei},
    title={Derived categories of twisted sheaves on elliptic threefolds},
    journal={J. Reine Angew. Math.},
    volume={544},
    date={2002},
    pages={161--179},
    issn={0075-4102},
}

%
%
%


\bib{dejong}{article}{
    author={de Jong, Aise Johan},
    title={A result of Gabber},
    eprint={http://www.math.columbia.edu/~dejong/},
}







 





 
%


%
%



%





%
%











 

\bib{lieblich-moduli}{article}{
    author={Lieblich, Max},
    title={Moduli of twisted sheaves},
    journal={Duke Math. J.},
    volume={138},
    date={2007},
    number={1},
    pages={23--118},
    issn={0012-7094},
}

%












\bib{mukai}{article}{
    author={Mukai, Shigeru},
    title={Duality between $D(X)$ and $D(\hat X)$ with its application to
    Picard sheaves},
    journal={Nagoya Math. J.},
    volume={81},
    date={1981},
    pages={153--175},
    issn={0027-7630},
}

%





%
%


\bib{rickard-morita}{article}{
    author={Rickard, Jeremy},
    title={Morita theory for derived categories},
    journal={J. London Math. Soc. (2)},
    volume={39},
    date={1989},
    number={3},
    pages={436--456},
    issn={0024-6107},
}

\bib{rickard-derived}{article}{
    author={Rickard, Jeremy},
    title={Derived equivalences as derived functors},
    journal={J. London Math. Soc. (2)},
    volume={43},
    date={1991},
    number={1},
    pages={37--48},
    issn={0024-6107},
}


 
%
%

 




 
 
 
\bib{thomason-trobaugh}{article}{
    author={Thomason, R. W.},
    author={Trobaugh, Thomas},
    title={Higher  algebraic   $K$-theory   of   schemes   and   of   derived   categories},
    conference={
    title={The Grothendieck Festschrift, Vol.\ III},
    },
    book={
        series={Progr. Math.},
        volume={88},
        publisher={Birkh\"auser Boston},
        place={Boston, MA},
    },
    date={1990},
    pages={247--435},
}


   

\bib{toen-derived}{article}{
    author = {To{\"e}n, Bertrand},
    title = {Derived Azumaya algebras and generators for twisted derived categories},
    journal = {Invent. Math.},
    year = {2012},
    volume = {189},
    number = {3},
    pages = {581--652},
}









\end{biblist}
\end{bibdiv}
\end{document}